\documentclass[11pt]{article}
\usepackage{graphicx}
\usepackage{amsmath}
 \vspace{0.6cm}

\newtheorem{lemma}{Lemma}[section]

\newtheorem{them}[lemma]{Theorem}

\newtheorem{cor}{Corollary}

\begin{document}

\title{\textbf{
Generalized Spectral Characterization of Graphs:
Revisited}\footnote{This work is supported by National Natural
Science Foundation of China 11071191.}}
\author{
\small Wei Wang \footnote{ E-mail
address: wang$\_$weiw@163.com. }
\\  \small School of Mathematics and Statistics, Xi'an Jiaotong University, Xi'an, 710049, P.R. China}
\date{}
\maketitle


\abstract \vspace{0.5cm} \noindent ~~~~A graph $G$ is said to be
\textit{determined by its generalized spectrum} (DGS for short) if
for any graph $H$, $H$ and $G$ are cospectral with cospectral
complements implies that $H$ is isomorphic to $G$. In
\cite{WX,WX1}, Wang and Xu gave some methods for determining
whether a family of graphs are DGS. In this paper, we shall review
some of the old results and present
some new ones along this line of research. 
\newpage
\section{Introduction}

Throughout the paper, we are only concerned with simple graphs.
Let $G$ be a simple graph with (0,1)-adjacency matrix $A(G)$.
\textit{The spectrum} of $G$ consists of all the eigenvalues
(together with their multiplicities) of the matrix $A(G)$. The
spectrum of $G$ together with that of its complement will be
referred to as \textit{the generalized spectrum} of $G$ in the
paper. For some terms and terminologies on graph spectra, see
\cite{CK}.

A graph $G$ is said to be \textit{determined by its spectrum} (DS
for short), if any graph having the same spectrum as $G$ is
necessarily isomorphic to $G$ (of course, the spectrum concerned
should be specified).

The spectrum of a graph encodes useful combinatorial information
about the given graph, and the relationship between the structural
properties of a graph and its spectrum has been studied
extensively over many years. A fundamental question in Spectral
Graph Theory is: `` Which graphs are DS?"  The problem dates back
to more than 50 years ago and originates from Chemistry, which has
received a lot of attention from researchers in recent years. It
turns out that, however, determining what kinds of graphs are DS
is generally a very hard problem.  For the background and some
known results about this problem, we refer the reader to
\cite{DH,DH1} and the references therein.

In \cite{WX,WX1}, Wang and Xu gave a method for determining
whether a graph $G$ is determined by its generalized spectrum (DGS
for short), by using some arithmetic properties of the walk-matrix
associated with the given graph. In this paper, we review some of
the previous results and further present some new results along
this line of research, which significantly improves the results in
\cite{WX,WX1}.
 The new ingredient of the paper is the discovery that whether the determinant of the walk-matrix is
square-free (for odd primes) is closely related to whether $G$ is
DGS.

The paper is organized as follows: In the next section, we review
some previous results that will be needed in the sequel. In
Section 3, we give a simple criterion for excluding odd primes.
The case $p=2$ is discussed in Section 4. Conclusions and open
problems are given in Section 5.


\section{Preliminaries}

For convenience of the reader, in this section, we will briefly
review some known results from \cite{WX,WX1}.

Let $W=[e,Ae,\cdots,A^{n-1}e]$ ($e$ is the all-one vector) be the
\textit{walk-matrix} of a graph $G$. Then the $(i,j)$-th entry of
$W$ is the number of walks of $G$ starting from vertex $i$ with
length $j-1$. A graph $G$ is called \textit{controllable graph} if
$W$ is non-singular (see also \cite{God}). It turns out that the
arithmetic properties of $det(W)$ is closely related to wether $G$
is DGS or not, as we shall see later. Denote by $\mathcal{G}_n$
the set of all controllable graphs on $n$ vertices. The following
theorem lies at the heart of our discussions.

\begin{them} [Wang and Xu \cite{WX}] Let $G\in{\mathcal{G}_n}$. Then there
exists a graph $H$ that is cospectral with $G$ w.r.t. the
generalized spectrum if and only if there exists a rational
orthogonal matrix $Q$ such that $Q^TA(G)Q=A(H)$ and $Qe=e$.
\end{them}

 Define
$$\mathcal{Q}_G=\left\{ \begin{array}{rrr}Q~is ~a~rational~~~&\vline&Qe=e,~Q^TA(G)Q ~is~a~ symmetric~\\
   ~orthogonal ~matrix&\vline&   (0,1)-matrix~ with~ zero~ diagonal
  \end{array}\right\},$$ where $e$ is the all-one vector. The following theorem follows easily from
  Theorem 1.1.
\begin{them}[Wang and Xu \cite{WX}] Let $G\in\mathcal{G}_n$. Then $G$ is DS w.r.t. the generalized spectrum iff the set $\mathcal{Q}_G$ contains only permutation matrices.
\end{them}

By the theorem above, in order to determine whether a given graph
$G\in\mathcal{G}_n$ is DS or not w.r.t. the generalized spectrum,
we have to determine those $Q$ in $\mathcal{Q}_G$ explicitly. At
first glance, this seems to be as difficult as the original
problem. However, we manage to do so by introducing the following
useful notion.

The \textit{level} of a rational orthogonal matrix $Q$ with $Qe=e$
is the smallest positive integer $\ell$ such that $\ell Q$ is an
integral matrix. Clearly, $\ell$ is the least common denominator
of all the entries of the matrix $Q$. If $\ell=1$, then $Q$ is a
permutation matrix.

Determining $\mathcal{Q}_G$ for all $G\in{\mathcal{G}_n}$ seems
too ambitious. Next, we shall only consider those controllable
graphs $G$ such that the level of those $Q\in{\mathcal{Q}_G}$
equals either 1 or 2.

 To illustrate the methods in \cite{WX,WX1}, first we give the relationships
between the values of $\ell$ for matrices $Q\in \mathcal{Q}_G$ and
properties of the walk-matrix $W$ of $G$. Recall that an $n\times
n$ matrix $U$ with integer entries is called \textit{unimodular}
if $det(U)=\pm 1 $. \textit{The Smith Normal Form} (SNF in short)
of an integral matrix $M$ is of the form
$diag(d_1,d_2,\cdots,d_n)$, where $d_i$ is the $i^{th}$
\textit{elementary divisor }of the matrix $M$ and
$d_i|d_{i+1}~(i=1,2,\cdots,n-1)$ hold. It is well known that for
every integral matrix $M$ with full rank, there exist unimodular
matrices $U$ and $V$ such that $
M=USV=Udiag(d_1,d_2,\cdots,d_n)V$, where $S$ is the SNF of the
matrix $M$. For a graph $G\in{\mathcal{G}_n}$, it is not difficult
to show that $d_n$ is the smallest positive integer $\ell$ such
that $\ell W^{-1}$ is an integral matrix.

\begin{them}[Wang and Xu \cite{WX}, Exclusion Criterion]
Let $W$ be the walk-matrix of a graph $G\in{\mathcal{G}_n}$, and
$Q\in{\mathcal{Q}_G}$ with level $\ell$. Then we have:

\noindent(a) $\ell|d_n$, where $d_n$ is the $n^{th}$ elementary
divisor of the SNF of $W$.

\noindent (b) Let $p$ be any prime factor of $d_n$. If $p|\ell$,
then the following system of congruence equations must have a
non-trivial solution ($x\not\equiv 0~mod~p$ ).
\begin{equation}\label{A1}
W^Tx\equiv 0,~x^Tx\equiv 0~(mod~ p).
\end{equation}
\end{them}

Theorem 2.3 (a) shows that $\ell$ is a divisor of $d_n$, and hence
all possible values of $\ell$ are finite for a given graph in
$\mathcal{G}_n$ and can be effectively computed through
calculating the SNF of $W$. While (b) shows that not all of the
divisors of $d_n$ can be a divisor of $\ell$; let $p$ be any prime
factor of $d_n(G)$ and if (\ref{A1}) has no non-trivial solution,
then $p$ must not be a prime factor of $\ell$, and it can be
excluded from further consideration. Using this way, it can be
expected that in most cases, many possibilities of the values of
the divisors of $d_n$ can be excluded.

Now we show how to check whether Eq. (\ref{A1}) has only trivial
solutions. As an illustration, we shall restrict ourselves to the
simplest case.

For convenience, we work with the finite field $\textbf{F}_p$ in
what follows. Suppose that $rank_p(W)=n-1$, where $rank_p(W)$ is
the rank of $W$ over the finite field \textbf{$\textbf{F}_p$}.
Consider the first equation of Eq. (\ref{A1}) as a system of
linear equations over $\textbf{F}_p$, then the set of solutions to
the first equation of (\ref{A1}) forms a one-dimensional subspace
of $\textbf{F}_p^n$. We can write $x=k\xi$, for some $0\neq
\xi\in{\textbf{F}_p^n}$ and $k=0,\cdots,p-1$. So Eq. (\ref{A1})
has only trivial solution iff
\begin{equation}\label{A2}
\xi^T\xi\neq 0~ in~ \textbf{F}_p.
\end{equation}

Let us give two examples which are taken from \cite{WX1}.

 Let $G_1$ and $G_2$ be two graphs with the adjacency matrices being given
as follows. It can easily be computed that $d_{12}(G_1)=2\cdot
17\cdot67\cdot8054231$, and $\xi^T\xi=12,25$ and $1492735$ for
each prime $p=17,67$ and $8054231$ respectively, where $\xi$ is
defined as above. Thus, all the prime factors of $d_{12}(G_1)$ can
be excluded except for $p=2$. It can be computed that
$d_{13}(G_2)=2\cdot3^2\cdot5\cdot197\cdot263\cdot5821$, and
$\xi^T\xi=1,0,139,101$ and $4298$ for each prime $p=3,5,197,263$
and $5821$. So all the prime factors of $d_{13}(G_2)$ can be
excluded except for $p=2,5$.

$$\tiny{\left[\begin{array}{cccccccccccc}0 & 1 & 1 & 0 & 0 & 1 & 0 & 0 & 1 & 0 & 1 & 1 \cr 1
& 0 & 0 & 0 & 0 & 0 & \ 0 & 0 & 1 & 1 & 0 & 1 \cr 1 & 0 & 0 & 0 &
0 & 1 & 0 & 0 & 1 & 0 & 1 & 0 \cr \ 0 & 0 & 0 & 0 & 1 & 0 & 0 & 1
& 0 & 1 & 0 & 0 \cr 0 & 0 & 0 & 1 & 0 & 1 & 0 \ & 0 & 0 & 0 & 0 &
1 \cr 1 & 0 & 1 & 0 & 1 & 0 & 1 & 1 & 0 & 1 & 1 & 0 \cr 0 \ & 0 &
0 & 0 & 0 & 1 & 0 & 0 & 1 & 1 & 1 & 0 \cr 0 & 0 & 0 & 1 & 0 & 1 &
0 & \ 0 & 1 & 0 & 1 & 0 \cr 1 & 1 & 1 & 0 & 0 & 0 & 1 & 1 & 0 & 0
& 0 & 0 \cr 0 & \ 1 & 0 & 1 & 0 & 1 & 1 & 0 & 0 & 0 & 0 & 1 \cr 1
& 0 & 1 & 0 & 0 & 1 & 1 & 1 \ & 0 & 0 & 0 & 1 \cr 1 & 1 & 0 & 0 &
1 & 0 & 0 & 0 & 0 & 1 & 1 & 0 \cr
\end{array}\right]~~~~\left[\begin{array}{ccccccccccccc}0 & 1 & 0 & 0 & 1 & 1 & 0 &
0 & 1 & 1 & 0 & 0 & 1 \cr 1 & 0 & 1 & 1 & 0 & \ 1 & 1 & 0 & 0 & 0
& 0 & 0 & 0 \cr 0 & 1 & 0 & 1 & 0 & 0 & 1 & 1 & 1 & 0 & 0 \ & 1 &
1 \cr 0 & 1 & 1 & 0 & 1 & 0 & 1 & 0 & 0 & 0 & 1 & 1 & 1 \cr 1 & 0
& 0 \ & 1 & 0 & 0 & 0 & 1 & 1 & 1 & 0 & 0 & 1 \cr 1 & 1 & 0 & 0 &
0 & 0 & 1 & 1 & \ 0 & 0 & 0 & 1 & 1 \cr 0 & 1 & 1 & 1 & 0 & 1 & 0
& 0 & 0 & 0 & 1 & 1 & 0 \cr \ 0 & 0 & 1 & 0 & 1 & 1 & 0 & 0 & 0 &
1 & 1 & 0 & 1 \cr 1 & 0 & 1 & 0 & 1 & 0 \ & 0 & 0 & 0 & 1 & 0 & 0
& 0 \cr 1 & 0 & 0 & 0 & 1 & 0 & 0 & 1 & 1 & 0 & 1 & \ 1 & 0 \cr 0
& 0 & 0 & 1 & 0 & 0 & 1 & 1 & 0 & 1 & 0 & 0 & 1 \cr 0 & 0 & 1 & \
1 & 0 & 1 & 1 & 0 & 0 & 1 & 0 & 0 & 1 \cr 1 & 0 & 1 & 1 & 1 & 1 &
0 & 1 & 0 \ & 0 & 1 & 1 & 0 \cr
  \end{array}\right]}.$$

 Nevertheless, it is not difficult to show that $p=2$ is always a
prime factor of $d_n$ and cannot be excluded invariably. In
\cite{WX1}, some further exclusion criterions are proposed to
eliminate the possibility of $p=2$. It can be show that $p=2$ can
be excluded for both graphs $G_1$ and $G_2$, by using the methods
in \cite{WX1}. Therefore $G_1$ is DGS. However, we do not know
wether $G_2$ is DGS or not since $p=5$ cannot be excluded using
the existing method.

In the next section, we shall present a simple criterion for
excluding primes $p>2$.

\section{A simple exclusion criterion for $p>2$}

In this section, we give a simple criterion for excluding primes
$p>2$, in terms of wether the exponent of $p$ in $det(W)$ is
larger than one. The main result of this section is the following
\begin{them}Let $G\in{\mathcal{G}_n}$, $Q\in{\mathcal{Q}_G}$ with
level $\ell$, and $p$ an odd prime. If $p|det(W)$ and
$p^2\not|det(W)$, then $p$ cannot be a divisor of $\ell$.
\end{them}

Before presenting the proof of above theorem, we need several
lemmas below. Note that the assumption that $p|det(W)$ and
$p^2\not|det(W)$ imply that $rank_p(W)=n-1$. This fact will be
used frequently in the sequel.

\begin{lemma}\label{LL1}Let $G\in{\mathcal{G}_n}$,
$Q\in{\mathcal{Q}_G}$ with level $\ell$. Let $p$ be an odd prime
divisor of $\ell$. Assume that $rank_p(W)=n-1$. Then we must have
$rank_p(\ell Q)=1$, and the following congruence equation has a
solution $z$:
\begin{equation}\label{M1}
Az\equiv \lambda_0z,~ e^Tz\equiv 0, z^Tz\equiv 0,~z\not\equiv
0~(mod~p)
\end{equation}
for some integer $\lambda_0$.
\end{lemma}
\textit{Proof.} The lemma follows immediately from the proof of
the next lemma.\hfill$\sharp$\\

\begin{lemma}\label{LL1}Let $G\in{\mathcal{G}_n}$,
$Q\in{\mathcal{Q}_G}$ with level $\ell$. Let $p$ be an odd prime
divisor of $\ell$. Assume that $rank_p(W)=n-1$ and $rank_p(\ell
Q)=1$, and the following congruent equation has a solution $z$:
\begin{equation}\label{M1}
W^Tz\equiv
0,~z^Tz\equiv 0,~ z\not\equiv 0~(mod ~p),
\end{equation}
Then $z^TAz\equiv \lambda_0 z^Tz ~(mod~p^2)$ holds, where
$\lambda_0$ is an integer such that $Az\equiv \lambda_0z~(mod~p)$
holds .
\end{lemma}
\textit{Proof.} First we claim that there exists a column $u$ of
the integer matrix $\ell Q$ and an integer vector $\beta$ such
that
\begin{equation}\label{E1}u=z+p\beta;
\end{equation}
\begin{equation}\label{E2}
u^TAu\equiv 0~(mod~p^2);
\end{equation}
\begin{equation}\label{E3}
u^Tu\equiv 0~(mod~p^2).
\end{equation}
In fact, it is easy to see that there exists a column $u$ of $\ell
Q$
 such that $u\not\equiv 0~(mod ~p)$.
With such a $u$, we have $W^Tu\equiv 0~(mod~p), u^Tu=\ell^2\equiv
0~(mod~p^2),$ and $u^TAu=0$. So $u$ is a solution of Eq.
(\ref{M1}), and Eq. (\ref{E1}) holds for some integer $\beta$.

 By Eq. (\ref{E3}) we have $$(z+p\beta)^T(z+p\beta)\equiv
z^Tz+2pz^T\beta\equiv 0~(mod~p^2).$$

Since $Q\in{\mathcal{Q}_G}$, we get $Q^TAQ=B$, where $B$ is the
adjacency matrix of some graph $H$. By $AQ=QB$ we get
$$Au_i=\sum_{k=1}^{n}b_{ik}u_k,i=1,2,\cdots,n,$$
where $u_i$ is the $i$-th column of $\ell Q$. Note that
$rank_p(\ell Q)=1$. Taking $mod ~p$ on both sides of the equation
above that contains $u$ on the right side generates $Az\equiv
\lambda_0 z~(mod~p)$, for some integer $\lambda_0$.

 Let $Az=\lambda_0z+p\gamma$, where $\gamma$ is an integer vector.
Then it follows from Eq. (\ref{E1}) and (\ref{E2}) that
\begin{eqnarray*}(z+p\beta)^TA(z+p\beta)&\equiv&z^TAz+2pz^TA\beta\\
&=&z^T(\lambda_0z+p\gamma)+2p(\lambda_0z+p\gamma)^T\beta \\
&\equiv &\lambda_0(z^Tz+2pz^T\beta)+pz^T\gamma\\
&\equiv& pz^T\gamma\\
&\equiv & 0~(mod~p^2)
\end{eqnarray*}
Thus we have $z^T(Az-\lambda_0z)=pz^T\gamma\equiv 0~(mod~p^2)$.
This completes the
proof. \hfill$\sharp$\\

\begin{lemma}\label{LE} Let $M=Udiag(d_1,d_2,\cdots,d_n)V=USV$, where $S$ is the SNF of
$M$, $U$ and $V$ are unimodular matrices and $d_i|d_{i+1}$ for
$i=1,2,\cdots,n-1$.  Then the system of congruence equations
$Mx\equiv 0~(mod~p^2)$ has a solution $x\not\equiv 0~(mod~p)$ if
and only if $p^2|d_n$.
\end{lemma}
\textit{Proof.} The equation $Mx\equiv 0~(mod~p^2)$ is equivalent
to $diag(d_1,d_2,\cdots,d_n)Vx\equiv 0~(mod~p^2)$. Let $Vx=y$.
Consider $diag(d_1,d_2,\cdots,d_n)y\equiv 0~(mod~p^2)$. If
$p^2|d_n$, let $y=(0,0,\cdots,0,1)^T$, then $x=V^{-1}y\not\equiv
0~(mod~p)$ is a required solution to the original congruence
equation. On the other hand, it is easy to see if $p^2\not| d_n$,
then the equation has no solution $x$ with $x\not\equiv
0~(mod~p)$.   \hfill$\sharp$\\

As a simple consequence of the above lemma, we have

\begin{cor}\label{LM1}Suppose that $rank_p(W)=n-1$, and $W^Tz\equiv 0,z\not\equiv
0~(mod~p)$. If there exists an integer vector $x$ such that
$W^Tx\equiv \frac{W^Tz}{p}~(mod~p)$, then $p^2|det(W)$.
\end{cor}

\begin{lemma}\label{P} If $rank_p(W)=n-1$, then we always have $rank_p(A-\lambda_0 I)\geq
n-2$.
\end{lemma}
\textit{Proof.} For contrary, suppose that there exist three
vectors $z, u$ and $v$ which are linearly independent over
$\textbf{F}_p$ such that $(A-\lambda_0 I)z=0,(A-\lambda_0 I)u=0$
and $(A-\lambda_0 I)v=0$, where we assume without loss of
generality that $e^Tz=0$, $e^Tu\neq 0$ and $e^Tv\neq 0$. Then we
can choose integers $k$ and $l$ with $ke^Tu+le^Tv=0$, over
$\textbf{F}_p$. Let $w=ku+lv$. Then $e^TA^iw=0$ and hence $W^Tw=0$
and $W^Tz=0$, which
implies that $rank_p(W)\leq n-2$, which contradicts the assumption that $rank_p(W)=n-1$. \hfill$\sharp$\\

It follows from Lemma \ref{P} that $rank_p(A-\lambda_0 I)=n-1$ or
$n-2$. Next, we shall distinguish this two cases in the following
lemmas.

\begin{lemma}If $rank_p(A-\lambda_0I)=n-1$, then
$p^2|det(W)$.
\end{lemma}
\textit{ Proof.} Let $z$ be an integral vector with $W^Tz\equiv
0~(mod~p)$. We prove the lemma by showing that the following
congruence equation
 always has a solution $x$.
\begin{equation}
 W^Tx\equiv\frac{W^Tz}{p}~(mod~p).
\end{equation}

Note that $z^Te=0$ and $z^T(A-\lambda_0I)=0$, over $\textbf{F}_p$.
It follows that the all-one vector $e$ can be written as the
linear combinations of the columns of $A-\lambda_0I$, i.e., there
exists a column vector $u$ such that
 \begin{equation}\label{Es1}
 e=(A-\lambda_0I)u,~\mbox{over}~\textbf{F}_p
\end{equation}
It follows from Eq. (\ref{Es1}) that there exists an integral
vector $\beta$ such that

 \begin{equation}\label{Eq.2} e=(A-\lambda_0I)u+p\beta.
 \end{equation}
Thus, we have
 \begin{eqnarray*}
 W&=&[e,Ae,\cdots,A^{n-1}e]\\
&=&[(A-\lambda_0I)u+p\beta,A((A-\lambda_0I)u+p\beta),\cdots,A^{n-1}((A-\lambda_0I)u+p\beta)]\\
&=&(A-\lambda_0I)[u,Au,\cdots,A^{n-1}u]+p[\beta,A\beta,\cdots,A^{n-1}\beta]\\
&=&(A-\lambda_0I)X+p[\beta,A\beta,\cdots,A^{n-1}\beta],
 \end{eqnarray*}
where $X:=[u,Au,\cdots,A^{n-1}u]$.

It follows that
 \begin{equation}\label{Es2}
 W^Tz=X^T(A-\lambda_0I)z+p[z^T\beta,z^TA\beta,\cdots,z^TA^{n-1}\beta]^T.
 \end{equation}

Since $W^Tz\equiv~ 0, (A-\lambda_0I)z\equiv ~0,~ A^iz\equiv
\lambda_0^{i}z~(i=0,1,\cdots,n-1)~(mod~p)$, by Eq. (\ref{Es2}) we
have
\begin{equation}\label{Eq4}
\frac{W^Tz}{p}\equiv
X^T\frac{(A-\lambda_0I)z}{p}+z^T\beta[1,\lambda_0,\cdots,\lambda_0^{n-1}]^T~(mod~p).
 \end{equation}

 Moreover, it follows from the fact that $rank_p(A-\lambda_0I)=n-1$, $z^T(A-\lambda_0I)=0$ and $z^Tz=0$, over
 $\textbf{F}_p$,
 that $z$ can be written as the linear combinations of the columns of $A-\lambda_0I$, i.e., there exists
 a vector $y$ such that $z=(A-\lambda_0I)y$.

 It is easy to show that $W^Ty\equiv
 e^Ty[1,\lambda_0,\cdots,\lambda_0^{n-1}]^T~ (mod~p).$ In fact,
 this follows from the following congruence equations:
 \begin{eqnarray*}
 z\equiv (A-\lambda_0I)y~(mod~p),~~~~~~~~~~~~\\
 e^TAy\equiv \lambda_0 e^Ty+e^Tz\equiv \lambda_0e^Ty~(mod~p),\\
 \cdots\cdots~~~~~~~~~~~~~~~~~~~~~~~~~~\\
 e^TA^{n-1}y\equiv \lambda_0^{n-1}e^Ty~(mod~p).
 \end{eqnarray*}

Now we show that $e^Ty\not\equiv 0~(mod~p)$. For otherwise, if
$e^Ty\equiv 0~(mod~p)$, then it follows that $W^Ty=0$ over
$\textbf{F}_p$. Note that $W^Tz=0$ over $\textbf{F}_p$. Moreover,
$y$ and $z$ are linearly independent. It follows that
$rank_p(W)\leq n-2$, which contradicts the fact that
$rank_p(W)=n-1$.

Thus, there exists an integer $k$ such that
\begin{equation}\label{EE1}
z^T\beta\equiv ke^Ty~(mod~p),
 \end{equation}

 Moreover, it follows from
the facts that $z^T\frac{(A-\lambda_0I)z}{p}\equiv 0,
~z^T(A-\lambda_0I)\equiv0,~(mod~p)$ and $rank_p(A-\lambda_0I)=n-1$
that the vector $\frac{(A-\lambda_0I)z}{p}$ can be written as the
linear combinations of the columns of $A-\lambda_0I$, i.e., there
exists a vector $v$ such that $$\frac{(A-\lambda_0I)z}{p}\equiv
(A-\lambda_0I)v.$$

 Note that $W^T\equiv X^T(A-\lambda_0I)~(mod~p)$.
Therefore, we have
 \begin{eqnarray*}
\frac{W^Tz}{p}&\equiv&X^T\frac{(A-\lambda_0I)z}{p}+kW^Ty~\\
&\equiv&W^Tv+kW^Ty\\
&\equiv& W^T(v+ky)~~~~~(mod~p).
 \end{eqnarray*}

By Cor. 1, the lemma follows.\hfill$\sharp$\\

\begin{lemma} \label{LQ} Let  $rank_p(W)=n-1$. Suppose that $rank_p(A-\lambda_0I)=n-2$. Then
$rank_p([A-\lambda_0I,z])=n-1$.
\end{lemma}
\textit{Proof.} Since $rank_p(A-\lambda_0I)=n-2$, there are two
vectors $z$ and $y$ which are linearly independent such that
$Az=\lambda_0z$ and $Ay=\lambda_0y$ with $e^Tz=0$, over
$\textbf{F}_p$.

Suppose the lemma does not hold. Then we have that $z$ can be
written as the linear combinations of the columns of
$A-\lambda_0I$. Thus, there exists a vector $x$ such that
$z=(A-\lambda_0I)x$, i.e.,
$$Ax=z+\lambda_0 x,$$
$$A^2x=Az+\lambda_0z+\lambda_0^2x,$$
$$\cdots\cdots$$
$$A^{n-1}x=A^{n-2}z+\lambda_0A^{n-3}z+\cdots+\lambda_0^{n-3}Az+\lambda_0^{n-2}z+\lambda_0^{n-1}x.$$
Now choose $k$ and $l$, not all zero, such that $e^Tw=0$, where
$w=kx+ly$.

Then, we have
$$e^TA^iw=ke^TA^ix+le^TA^iy=ke^T(A^{i}z+\lambda_0A^{i-1}z+\cdots+\lambda_0^{i-1}Az+\lambda_0^{i-1}z)+\lambda_0^{i}(ke^Tx+le^Ty)=0,$$
for $i=0,1,\cdots,n-1$, i.e., $W^Tw=0$.

Now we show that $x,y$ and $z$ are linearly independent. Suppose
$ax+by+cz=0$. Then left-multiplying both sides of the above
equality by $(A-\lambda_0I)$ gives $az=0$, which implies $a=0$. By
assumption that $y$ and $z$ are linearly independent, we have
$b=c=0$.

  Therefore, $z$ and $w$ are linearly independent.
Moreover, we have $W^Tz=0$ and $W^Tw=0$. This contradicts the fact that $rank_p(W)=n-1$. \hfill$\sharp$\\

\begin{lemma}Suppose that $rank_p(A-\lambda_0I)=n-2$. Then
$p^2|det(W)$
\end{lemma}
\textit{Proof.} Note that $rank_p(W)=n-1$ and
$rank_p(A-\lambda_0I)=n-2$. By Lemma \ref{LQ}, we get that $z$
cannot be expressed as the linear combinations of the column
vectors of $A-\lambda_0I$, over $\textbf{F}_p$, and hence
$rank_p([A-\lambda_0,z])=n-1$. Moreover, $z^Te=0$ and
$z^T[A-\lambda_0I,z]=0$, it follows that the all-one vector $e$
can be expressed as the linear combinations of the column vectors
of $A-\lambda_0I$ and $z$, i.e., there exist an vector $u$ and an
integer $k$ such that

$$e=(A-\lambda_0I)u+kz, ~over~\textbf{F}_p.$$  That is,
$$e=(A-\lambda_0I)u+kz+p\beta,~ over~ \textbf{Z}.$$

It follows that
 \begin{eqnarray*}
 Ae&=&A(A-\lambda_0I)u+kAz+pA\beta=(A-\lambda_0I)Au+kAz+pA\beta.\\
A^2e&=&A^2(A-\lambda_0I)u+kA^2z+pA^2\beta=(A-\lambda_0I)A^2u+kA^2z+pA^2\beta.\\
&&~~~~~~~~~~~~~~~~\cdots\cdots\cdots\cdots\cdots\\
A^{n-1}e&=&A^{n-1}(A-\lambda_0I)u+kA^{n-1}z+pA^{n-1}\beta=(A-\lambda_0I)A^{n-1}u+kA^{n-1}z+pA^{n-1}\beta.
 \end{eqnarray*}
Therefore,
\begin{eqnarray*}
 W&=&[e,Ae,\cdots,A^{n-1}e]\\
&=&(A-\lambda_0I)[u,Au,\cdots,A^{n-1}u]+k[z,Az,\cdots,A^{n-1}z]+p[\beta,A\beta,\cdots,A^{n-1}\beta]\\
&=&(A-\lambda_0I)X+k[z,Az,\cdots,A^{n-1}z]+p[\beta,A\beta,\cdots,A^{n-1}\beta],~over~\textbf{Z},
 \end{eqnarray*}
where $X=[u,Au,\cdots,A^{n-1}u]$. It follows that
 \begin{eqnarray*}
\frac{W^Tz}{p}&=&X^T\frac{(A-\lambda_0I)z}{p}+k[\frac{z^Tz}p,\frac{z^TAz}{p},\cdots,\frac{z^TA^{n-1}z}{p}]^T\\
&+&[\beta^Tz,\cdots,\beta^TA^{n-1}z]^T~~(over~\textbf{Z})\\
&\equiv&X^T\frac{(A-\lambda_0I)z}{p}+k\frac{z^Tz}{p}[1,\lambda_0,\cdots,\lambda_0^{n-1}]^T+\beta^Tz[1,\lambda_0,\cdots,\lambda_0^{n-1}]^T~(mod~p)
 \end{eqnarray*}
The congruence equation follows from the facts that
$\frac{z^TA^iz}{p}-\frac{\lambda_0^{i}z^Tz}{p}~\equiv0~\mbox{and}~A^iz\equiv
\lambda_0^iz~(mod~p)$.

Moreover, $z^T\frac{(A-\lambda_0I)z}{p}\equiv
0~\mbox{and}~z^T[A-\lambda_0I,z]\equiv 0~(mod~p)$. It follows that
there exist a vector $\alpha$ and an integer $m$ such that
\begin{equation}
\frac{(A-\lambda_0I)z}{p}\equiv (A-\lambda_0I)\alpha+mz~(mod~p)
\end{equation}
 \begin{eqnarray*}
X^T\frac{(A-\lambda_0I)z}{p}&\equiv&X^T(A-\lambda_0I)\alpha+mX^Tz\\
&\equiv& W^T\alpha-kz^T\alpha[1,\lambda_0,\cdots,\lambda_0^{n-1}]^T+mX^Tz\\
&\equiv&
W^T\alpha+(mu^Tz-kz^T\alpha)[1,\lambda_0,\cdots,\lambda_0^{n-1}]^T~~(mod~p)
\end{eqnarray*}
Thus
\begin{equation}
\frac{W^Tz}{p}\equiv W^T\alpha+(k\frac{z^Tz}{p}+\beta^Tz+
mu^Tz-kz^T\alpha)[1,\lambda_0,\cdots,\lambda_0^{n-1}]^T~~(mod~p)
\end{equation}
Let $y$ be a vector with $(A-\lambda_0)y=0$ that is linearly
independent with $z$. Then we must have $e^Ty\not\equiv~0(mod~p)$.
For otherwise, if $e^Ty\equiv~0(mod~p)$, then it follows
$W^Ty\equiv 0$. Note $W^Tz=0$, $W^Ty=0$ and $y$ and $z$ are
linearly independent, over $\textbf{F}_p$. This contradicts fact
that $rank_p(W)=n-1$.

 It follows that there exists an integer $l$ such that $$k\frac{z^Tz}{p}+\beta^Tz+
mu^Tz-kz^T\alpha\equiv le^Ty~(mod~p).$$ Thus, we have
\begin{equation} \frac{W^Tz}{p}\equiv
W^T\alpha+le^Ty[1,\lambda_0,\cdots,\lambda_0^{n-1}]^T\equiv
W^T\alpha+lW^Ty\equiv W^T(\alpha+ly) ~~(mod~p)
\end{equation}
By Cor. 1, the lemma follows.\hfill$\sharp$\\

\textit{\textbf{Proof of Theorem 3.1}}: Combining Lemmas 3.2-3.8,
Theorem 3.1 follows immediately. \hfill$\sharp$\\

Let us give a few remarks to end this section.

i) Result in Theorem 3.1 is best possible in the sense that if
$p>2$ has exponent larger than one, then Theorem 3.1 may not be
true. The following is a counterexample.

Let the adjacency matrix of graph $G$ be given as below. It can
easily be computed that
$$det(W)=2^6\times3^2\times157\times1361\times2237.$$
The exponent of $p=3$ in the standard prime decomposition $det(W)$
is equal to 2, and $p=3$ cannot be excluded. Actually, let $Q$ be
a rational orthogonal matrix given as below. Then
$Q\in{\mathcal{Q}_G}$ with level $\ell=3$, since it can be easily
verified that $Q^TAQ$ is a $(0,1)$-matrix.

$$\tiny{A=\left[\begin{array}{cccccccccccc} 0 & 0 & 0 & 0 & 0 & 1 & 0 & 1 & 0 & 0 & 1 & 0 \cr 0 & 0
& 1 & 0 & 1 & 0 & \ 1 & 0 & 1 & 0 & 1 & 1 \cr 0 & 1 & 0 & 1 & 1 &
0 & 1 & 1 & 0 & 0 & 1 & 0 \cr \ 0 & 0 & 1 & 0 & 1 & 0 & 1 & 0 & 0
& 0 & 1 & 1 \cr 0 & 1 & 1 & 1 & 0 & 1 & 0 \ & 1 & 0 & 0 & 0 & 1
\cr 1 & 0 & 0 & 0 & 1 & 0 & 1 & 1 & 1 & 0 & 1 & 0 \cr 0 \ & 1 & 1
& 1 & 0 & 1 & 0 & 1 & 1 & 0 & 1 & 1 \cr 1 & 0 & 1 & 0 & 1 & 1 & 1
& \ 0 & 1 & 1 & 0 & 0 \cr 0 & 1 & 0 & 0 & 0 & 1 & 1 & 1 & 0 & 0 &
0 & 1 \cr 0 & \ 0 & 0 & 0 & 0 & 0 & 0 & 1 & 0 & 0 & 1 & 1 \cr 1 &
1 & 1 & 1 & 0 & 1 & 1 & 0 \ & 0 & 1 & 0 & 1 \cr 0 & 1 & 0 & 1 & 1
& 0 & 1 & 0 & 1 & 1 & 1 & 0 \cr
\end{array}\right]},$$
$$\tiny{Q=\frac{1}{3}\left[\begin{array}{rrrrrrrrrrrr} 0 & 0 & 0 & 0 & 0 & 0 & 3 & 0 & 0 & 0 & 0 & 0 \cr 2 & -1 & \
-1 & 1 & 1 & 1 & 0 & 0 & 0 & 0 & 0 & 0 \cr 0 & 0 & 0 & 0 & 0 & 0 &
0 & 3 & \ 0 & 0 & 0 & 0 \cr 0 & 0 & 0 & 0 & 0 & 0 & 0 & 0 & 3 & 0
& 0 & 0 \cr 0 & 0 & \ 0 & 0 & 0 & 0 & 0 & 0 & 0 & 3 & 0 & 0 \cr 1
& 1 & 1 & 2 & -1 & \ -1 & 0 & 0 & 0 & 0 & 0 & 0 \cr -1 & 2 & \ -1
& 1 & 1 & 1 & 0 & 0 & 0 & 0 & 0 & 0 \cr -1 & \ -1 & 2 & 1 & 1 & 1
& 0 & 0 & 0 & 0 & 0 & 0 \cr 1 & 1 & 1 & -1 & 2 & \ -1 & 0 & 0 & 0
& 0 & 0 & 0 \cr 0 & 0 & 0 & 0 & 0 & 0 & 0 & 0 & 0 & 0 & 3 & \ 0
\cr 0 & 0 & 0 & 0 & 0 & 0 & 0 & 0 & 0 & 0 & 0 & 3 \cr 1 & 1 & 1 &
-1 & \ -1 & 2 & 0 & 0 & 0 & 0 & 0 & 0 \cr
\end{array}\right]}.$$\\

ii) By Theorem 3.1, for graph $G_2$ in the previous example, $p=5$
can also be excluded since the $5|det(W)$ and $5^2\not|det(W)$.
Thus, $G_2$ is also DGS.


\section{Some discussions on $p=2$}

As mentioned previously, the case $p=2$ is more involved to deal
with. Let us try to explain this through the following lemmas.

\begin{lemma}[c.f. Wang \cite{WX2}] $e^TA^ke$ is even for every positive integer $k$.
\end{lemma}
\textit{Proof.} Note that $$e^TA^ke=Tr(A^k)+\sum_{i\neq
j}A^k=Tr(A^k)+2\sum_{i< j}A^k\equiv Tr(A^k)~(mod~2).$$

$Tr(A^k)=Tr(AA^{k-1})=\sum_{i,j}a_{ij}b_{ij}=2\sum_{i<j}a_{ij}b_{ij}$,
where $B:=A^{k-1}$. Thus the lemma follows.\hfill$\sharp$\\

\begin{lemma}[c.f. Wang \cite{WX2}] $rank_2(W)\leq \lceil \frac{n}{2}\rceil.$
\end{lemma}
\textit{Proof.} Suppose $n$ is even. Then it follows from Lemma
4.1 that $W^TW=0$ over $\textbf{F}_2$.
$2rank_2(W)=rank_2(W^T)+rank_2(W)\leq n$. Thus we have
$rank_2(W)\leq n/2=\lceil \frac{n}{2}\rceil.$

If $n$ is odd. Let $\hat{W}$ be the matrix obtained from $W$ by
deleting the first column. Then $W^T\hat{W}=0$ over
$\textbf{F}_2$. Note $rank_2(W)+rank_2(\hat{W})\leq n$ and
$rank_2(\hat{W})\geq rank_2(W)-1$. It follows that $rank_2(W)\leq
(n+1)/2=\lceil \frac{n}{2}\rceil$.   \hfill$\sharp$\\

By Lemma 4.2, the system of linear equations in Eq. (1) has a set
of solutions with dimension at least $\lfloor n/2\rfloor $, and it
not difficult to show that it is always possible to choose some of
the solutions to meet the second requirement in Eq. (1).

Moreover, by Lemma 4.2, the following corollary follows
immediately.

\begin{cor} Let $det(W)=\epsilon 2^\alpha p_1^{\alpha_1}\cdots
p_s^{\alpha_s}$ ($\epsilon=\pm 1$) be the standard decomposition
of $det(W)$ into prime factors . Then $\alpha\geq \lfloor
\frac{n}{2}\rfloor$.
\end{cor}

For any graph $G\in{\mathcal{G}_n}$, the number of $d_i$ which is
even in the SNF $S=diag(d_1,d_2,\cdots,d_n)$ of $W$ must be at
leat $\lfloor n/2\rfloor$. Next, we are interested in a specific
family of controllable graphs
$$\mathcal{F}_n=\{G\in{\mathcal{G}_n}|\frac{det(W)}{2^{\lfloor n/2\rfloor}}\mbox{~is
~square-free~and}~2^{\lfloor n/2\rfloor+1}\not|det(W)\}.$$

By Cor. 2, for every graph in $\mathcal{F}_n$, the SNF of $W$ must
be like $S=diag(1,\cdots,1,2,\cdots,2,2b)$, where $b$ is an odd
square-free integer and the number of $2$'s is exactly $\lfloor
n/2\rfloor$ in the diagonal of $W$.

Let $G\in{\mathcal{F}_n}$. Let $Q\in{\mathcal{Q}_G}$ with level
$\ell$ and $p$ be any prime divisor of $\ell$. Then by Theorem 2.3
(a), we have $p|2b$. If $p>2$, then by Theorem 3.1, we have $p\not
|\ell$. Therefore, $\ell=1$ or $\ell=2$. Next, we present a simple
exclusion criterion for $\ell=2$, which significantly simplifies
the method in \cite{WX1}.

\begin{lemma}Let $G\in{\mathcal{G}_n}$. Let $Q\in{\mathcal{Q}_G}$ with level $\ell=2$. Then there exists a (0,1)-vector
$u$ with four non-zero entries `1' such that
 \begin{equation}u^TA^ku\equiv 0~(mod~4)~,k=1,2,\cdots,n-1.
\end{equation}
Moreover, $u$ satisfies $W^Tu\equiv 0,u\not\equiv 0~(mod~2)$
\end{lemma}
\textit{Proof. }$Q\in{\mathcal{Q}_G}$ implies that $Q^TAQ=B$,
where $B$ is a (0,1)-matrix. Let $\bar{u}$ be the $i$-th column of
$2Q$. It follows from $Q^TA^kQ=B^k$ that
$\bar{u}^TA^k\bar{u}=4(B^k)_{i,i}\equiv 0~(mod ~4)$. It follows
from the facts $\ell=2$ and $Qe=e$ that the four non-zero entries
of $\bar{u}$ are $1,1,1,$ and $-1$, respectively. Let
$u=\bar{u}+2e_j$ ($e_j$ denotes the $j$-th standard basis of
$\textbf{R}^n$) be a $(0,1)$-vector with four non-zero entries `1'
. Then
$$u^TA^ku=\bar{u}^TA^k\bar{u}+4\bar{u}^TA^ke_j+4e_j^TA^ke_j\equiv
0~(mod~4).$$

The last assertion follows from the fact that $Q^TA^kQ=B^k$ and
$Qe=e$ imply that $W^TQ$ is an integral matrix. Thus $W^Tu\equiv
0,u\not\equiv 0~(mod~2)$ holds.   \hfill$\sharp$\\

Lemma 4.3 gives a simple way to eliminate the possibility of
$\ell=2$. First, solve the system of linear equations $W^Tx=0$
with additional requirement that $x$ has four non-zero entries
$1$, over $\textbf{F}_2$, to get a solution set $S$. This can be
done through checking $n\choose 4$ possibilities. Then for each
solution $x$ check whether Eq. (17) holds. If all $x\in S$ does
not satisfy Eq. (17), then $\ell\neq 2$ and hence $\ell=1$, i.e.,
$G$ is DGS.

Let us give an example for illustration. Let $G=G_1$ be the first
graph given in Section 2. Clearly $G\in{\mathcal{F}_n}$. It can be
easily computed by Mathematica 5.0 that the corresponding solution
set is

$S=\{(0, 1, 0, 1, 0, 1, 0, 1, 0, 0, 0, 0)^T,(0, 0, 1, 0, 1, 0,
1, 0, 1, 0, 0, 0)^T,({1, 0, 0, 0, 1, 0, 1, 0, 0, 1, 0, 0})^T,\\
({1, 0, 1, 0, 0, 0, 0, 0, 1, 1, 0, 0})^T,
 ({1, 0, 1, 0, 0, 0, 1, 0, 0, 0, 1, 0})^T,({1, 0, 0, 0, 1, 0, 0,
0, 1, 0, 1, 0})^T,\\({0, 0, 1, 0, 1, 0, 0, 0, 0, 1, 1, 0})^T, ({0,
0, 0, 0, 0, 0, 1, 0, 1, 1, 1, 0})^T, ({1, 0, 1, 0, 1, 0, 0, 0, 0,
0, 0, 1})^T,\\({1, 0, 0, 0, 0, 0, 1, 0, 1, 0, 0, 1})^T, ({0, 0, 1,
0, 0, 0, 1, 0, 0, 1, 0, 1})^T,({0, 0, 0, 0, 1, 0, 0, 0, 1, 1, 0,
1})^T,\\({0, 0, 0, 0, 1, 0, 1, 0, 0, 0, 1, 1})^T,({0, 0, 1, 0, 0,
0, 0, 0, 1, 0, 1, 1})^T,({1, 0, 0, 0, 0, 0, 0, 0, 0, 1, 1,
1})^T\}$

However, non of $x\in{S}$ satisfies Eq. (17). Thus $G$ is DGS.

We remark, though Lemma 4.3 is a sufficient condition to exclude
the case $\ell=2$, our numerical experiments do suggest that it is
always necessary for graphs $G\in{\mathcal{F}_n}$.

\section{Concluding remarks and open problems}

We have reviewed some previous results on the topic of
characterizing a graph by both its spectrum and the spectrum of
its complement. Then we have presented a simple new exclusion
criterion for excluding odd primes. The case $p=2$ has also been
discussed.

As it turns out, the arithmetic properties of $det(W)$ is closely
related to whether a given controllable graphs is DGS. Actually,
we have the following\\

\textbf{Conjecture (Wang \cite{WX2})}: \textit{Every graph in
$\mathcal{F}_n$ is DGS.}\\

For a given graph $G\in{\mathcal{F}_n}$, $Q\in{\mathcal{Q}_G}$
with level $\ell$. We have shown that either $\ell=1$ or $\ell=2$.
However, some additional efforts have to be made to eliminate the
possibility of $\ell=2$.

 Finally, we remark that it can be shown
 (see \cite{WX3}) that almost every graphs in $\mathcal{F}_n$ is
 DGS. In view of the simple definition of $\mathcal{F}_n$, it
 suggests a possible way to show that DGS-graphs have positive density via
 proving $\mathcal{F}_n$ has positive density (numerical experiments show that $\mathcal{F}_n$ has density nearly 0.2).
This needs further investigations in the future.

%

\end{document}